\magnification=\magstep1
\tolerance=2000

\baselineskip=17pt

\font\tenopen=msbm10
\font\sevenopen=msbm7
\font\fiveopen=msbm5
\newfam\openfam
\def\openo{\fam\openfam\tenopen}
\textfont\openfam=\tenopen
\scriptfont\openfam=\sevenopen
\scriptscriptfont\openfam=\fiveopen

\font\tenujsym=msam10
\font\sevenujsym=msam7
\font\fiveujsym=msam5
\newfam\ujsymfam

\textfont\ujsymfam=\tenujsym
\scriptfont\ujsymfam=\sevenujsym
\scriptscriptfont\ujsymfam=\fiveujsym
\def\loc{{\rm loc}}
\def\e{\varepsilon} 

\def\R{{\openo R}}
\def\C{{\openo C}}
\def\N{{\openo N}}
\def\Z{{\openo Z}}

\def\kocka{\kern+0.9pt{\sqcup\kern-6.5pt\sqcap}\kern+0.9pt}

\def\zmatrix#1{\null\,\vcenter{\normalbaselines  
    \ialign{\hfil$##$\hfil&&\ \hfil$##$\hfil\crcr  
      \mathstrut\crcr\noalign{\kern-\baselineskip}%
      #1\crcr\mathstrut\crcr\noalign{\kern-\baselineskip}}}\,}%

\def\vonal{\noalign{\hrule}}

\def\zv{\vrule depth6pt height11pt }

\def\dddots{\mathinner{\mkern1mu\mkern2mu\raise1pt\hbox{.}
\mkern2mu\raise4pt\hbox{.}\raise7pt\vbox{\kern7pt\hbox{.}}
\mkern1mu}}%

\def\abra#1{\smallskip\centerline{\it Figure #1}\medbreak}

\centerline{\bf On the convergence of double integrals and} 

\centerline{\bf a generalized version of Fubini's theorem on successive integration} 

\bigskip

\centerline{FERENC M\'ORICZ}

\bigskip

\noindent Bolyai Institute, University of Szeged, Aradi v\'ertan\'uk tere 1, 

\noindent Szeged 6720, Hungary 

\noindent e-mail: moricz@math.u-szeged.hu

\vglue1cm

\noindent {\bf Abstract.} Let the function $f: \overline{\R}^2_+ \to \C$ be such that $f\in L^1_{\loc}  (\overline{\R}^2_+)$. We 
investigate the convergence behavior of the double integral 
$$\int^A_0 \int^B_0 f(u,v) du dv \quad {\rm as} \quad A,B \to \infty,\leqno(*)$$
where $A$ and $B$ tend to infinity independently of one another; 
while using two notions of convergence: that in Pringsheim's sense and that in the regular sense. 
Our main result is the following Theorem 3: If the double integral ($*$) converges in the 
regular sense, or briefly: converges regularly, then the finite limits 
$$\lim_{y\to \infty} \int^A_0 \Big(\int^y_0 f(u,v) dv\Big) du =: I_1 (A)$$
and
$$\lim_{x\to \infty} \int^B_0 \Big(\int^x_0 f(u,v) du) dv = : I_2 (B)$$
exist uniformly in $0<A, B <\infty$, respectively; and 
$$\lim_{A\to \infty} I_1(A) = \lim_{B\to \infty} I_2 (B) = 
\lim_{A, B \to \infty} \int^A_0 \int^B_0 f(u,v) du dv.$$
This can be considered as a generalized version of Fubini's theorem on successive integration when 
$f\in L^1_{\loc} (\overline{\R}^2_+)$, but 
$f\not\in L^1 (\overline{\R}^2_+)$. 

\bigskip


\bigskip

\noindent {\it Key words and phrases}: double series of complex numbers, double integrals  of 
locally integrable functions over $\overline{\R}^2_+$ in Lebesgue's sense, 
convergence in Pringsheim's sense, regular convergence, absolute convergence, a 
generalized version of Fubini's theorem on successive 
integration.

\vglue1cm 

\noindent {\bf 1. Background: Convergence of double series of complex numbers} 

We consider the double series 
$$\sum^\infty_{j=0} \ \sum^\infty_{k=0} a_{j,k}\leqno(1.1)$$
of complex numbers, in symbols: $\{a_{j,k} : (j,k) \in \N^2\} \subset \C$. As usual, the 
{\it rectangular partial sums} of (1.1) are defined by 
$$s_{m,n} : = \sum^m_{j=0} \ \sum^n_{k=0} a_{j,k}, \quad (m,n) \in \N^2.\leqno(1.2)$$

We recall (see [6] by A. Pringsheim) that the double series (1.1) is said to {\it converge} {\it in Pringsheim's sense} to the sum $s\in \C$, 
in symbols: 
$$\lim_{m,n \to \infty} \ s_{m,n} = s$$
if for every $\e > 0$ there exists 
$\kappa_1 = \kappa_1 (\e) \in \N$ such that 
$$|s_{m,n} - s| < \e \quad {\rm if} \quad \min \{m,n\} > \kappa_1. \leqno(1.3)$$
We note that A. Zygmund uses this convergence notion without the term ``in Pringsheim's sense"; see in 
[9, on p. 303, just after the formula (1.18), and he actually defines it for $d$-multiple series,  $d\ge 2$]. 

It has been observed by many authors that the convergence of a double series in Pringsheim's sense implies neither the boundedness of its terms 
$a_{j,k}$, nor the convergence of its so-called `row' and `column' series defined respectively by 
$$\sum^\infty_{k=0} a_{j,k}, \quad j\in \N; \quad {\rm and} \quad \sum^\infty_{j=0} a_{j,k}, \quad k\in \N.\leqno(1.4)$$
For the reader's convenience, we give the following Examples 1 and 2 of these phenomena. 

\bigskip 

\noindent {\bf Example 1.} Let the terms $a_{j,k}$ of (1.1) be given in Figure 1 below. Clearly, we have 
$$s_{m,n} = 0 \quad {\rm if} \quad \min \{m,n\} \ge 1.$$
Thus, in this case (1.1) converges to 0  in Pringsheim's sense, but its terms 
$a_{j,k}$ are not bounded. 
\bigskip

$$
\zmatrix{
\zv & \vdots & \zv & & \zv & & \zv  & & \zv & & \zv & & \zv & \dddots  \cr
\vonal
\zv & a_{5,0} & \zv & a_{5,1} & \zv & a_{5,2} & \zv & a_{5,3} & \zv & a_{5,4} & \zv & a_{5,5} & \zv  & \cr
\vonal
\zv & a_{4,0} & \zv & a_{4,1} & \zv & a_{4,2} & \zv &a_{4,3} & \zv & a_{4,4} & \zv & a_{4,5} & \zv  & \cr
\vonal
\zv & a_{3,0} & \zv & a_{3,1} & \zv & a_{3,2} & \zv & a_{3,3} & \zv & a_{3,4} & \zv & a_{3,5} & \zv  & \cr
\vonal
\zv & a_{2,0} & \zv & a_{2,1} & \zv & a_{2,2} & \zv & a_{2,3} & \zv & a_{2,4} & \zv & a_{2,5} & \zv  & \cr
\vonal
\zv & a_{1,0} & \zv & a_{1,1} & \zv & a_{1,2} & \zv & a_{1,3} & \zv & a_{1,4} & \zv & a_{1,5} & \zv  & \cr
\vonal
\zv & a_{0,0} & \zv & a_{0,1} & \zv & a_{0,2} & \zv & a_{0,3} & \zv & a_{0,4} & \zv & a_{0,5} & \zv & \cdots & \cr
\vonal
} \qquad\qquad
\zmatrix{
\zv & \vdots & \zv & & \zv & & \zv  & & \zv & & \zv & & \zv & \dddots  \cr
\vonal
\zv & -3 & \zv & 3 & \zv & 0 & \zv & 0 & \zv & 0 & \zv & 0 & \zv  & \cr
\vonal
\zv & 3 & \zv & -3 & \zv & 0 & \zv & 0 & \zv & 0 & \zv & 0 & \zv  & \cr
\vonal
\zv & -2 & \zv & 2 & \zv & 0 & \zv & 0 & \zv & 0 & \zv & 0 & \zv  & \cr
\vonal
\zv & 2 & \zv & -2 & \zv & 0 & \zv & 0 & \zv & 0 & \zv & 0 & \zv  & \cr
\vonal
\zv & -1 & \zv & 1 & \zv & -2 & \zv & 2 & \zv & -3 & \zv & 3 & \zv  & \cr
\vonal
\zv & 1 & \zv & -1 & \zv & 2 & \zv & -2 & \zv & 3 & \zv & -3 & \zv & \cdots & \cr
\vonal
}
$$
\abra {1}

\bigskip

\noindent {\bf Example 2.} Let the terms $a_{j,k}$ of (1.1) be given in Figure 2. In this case, each row 
series and each column series of (1.1) 
contains only two nonzero terms, but (1.1) fails to converge in Pringsheim's sense, since 
$$s_{m,m} = \cases{1 &if \quad $m$ \ {\rm is\ even},\cr
0 &if \quad $m$ \ \ {\rm is \ odd}, \ \ $m\in \N$. \cr}$$
\bigskip

$$
\zmatrix{
\zv & \vdots & \zv & & \zv & & \zv  & & \zv & & \zv & & \zv & \dddots  \cr
\vonal
\zv & 0 & \zv & 0 & \zv & 0 & \zv & 0 & \zv & -1 & \zv & 1 & \zv  & \cr
\vonal
\zv & 0 & \zv & 0 & \zv & 0 & \zv & 0 & \zv & 1 & \zv & -1 & \zv  & \cr
\vonal
\zv & 0 & \zv & 0 & \zv & -1 & \zv & 1 & \zv & 0 & \zv & 0 & \zv  & \cr
\vonal
\zv & 0 & \zv & 0 & \zv & 1 & \zv & -1 & \zv & 0 & \zv & 0 & \zv  & \cr
\vonal
\zv & -1 & \zv & 1 & \zv & 0 & \zv & 0 & \zv & 0 & \zv & 0 & \zv  & \cr
\vonal
\zv & 1 & \zv & -1 & \zv & 0 & \zv & 0 & \zv & 0 & \zv & 0 & \zv & \cdots & \cr
\vonal
}
$$
\abra {2}

\bigskip

The converse situation may also occur when the double series (1.1) converges in Prings\-heim's sense, while 
each row series and each column series of (1.1) diverges. See, e.g., the following 

\bigskip

\noindent {\bf Example 3.} Let the terms $a_{j,k}$ of (1.1) be given in Figure 3. In this case 
$$s_{m,n} = \cases{{2\over 2 + \min \{m,n\} } &if \quad {\rm both}\ $m$ \ {\rm and} \ $n$ \ {\rm are\ even},\cr
0 &otherwise.\cr}$$
This time  (1.1) converges to 0 in Pringsheim's sense, but each row series and each column series  diverges. 

\bigskip

$$
\zmatrix{
\zv & \vdots & \zv & & \zv & & \zv  & & \zv & & \zv & & \zv & \dddots  \cr
\vonal
\zv & -1 & \zv & 1 & \zv & -{1\over 2} & \zv & {1\over 2} & \zv & -{1\over 3} & \zv & {1\over 3} & \zv  & \cr
\vonal
\zv & 1 & \zv & -1 & \zv & {1\over 2} & \zv & -{1\over 2} & \zv & {1\over 3} & \zv & -{1\over 3} & \zv  & \cr
\vonal
\zv & -1 & \zv & 1 & \zv & -{1\over 2} & \zv & {1\over 2} & \zv & -{1\over 2} & \zv & {1\over 2} & \zv  & \cr
\vonal
\zv & 1 & \zv & -1 & \zv & {1\over 2} & \zv & -{1\over 2} & \zv & {1\over 2} & \zv & -{1\over 2} & \zv  & \cr
\vonal
\zv & -1 & \zv & 1 & \zv & -1 & \zv & 1 & \zv & -1 & \zv & 1 & \zv  & \cr
\vonal
\zv & 1 & \zv & -1 & \zv & 1 & \zv & -1 & \zv & 1 & \zv & -1 & \zv & \cdots & \cr
\vonal
}
$$
\abra {3}

\bigskip

We observe that if the double series (1.1) converges in Pringsheim's sense, then it follows from (1.2) and (1.3) that 
for every $\e>0$, 
$$\Big|\sum^M_{j=m} \ \sum^N_{k=n} a_{j,k}\Big| = |s_{M,N} - s_{m-1, N} - s_{M, n-1} + 
s_{m-1, n-1}| < 4\e\leqno(1.5)$$
$${\rm if} \quad \min\{m,n\} > \kappa_1(\e)\quad {\rm and} \quad M\ge m, \ N\ge n.$$
In contrast to single sequences, the fulfillment of (1.5) for every $\e>0$ does not imply the convergence 
of (1.1) in Pringsheim's sense. See, e.g., the following

\bigskip

\noindent {\bf Example 4.} Let the terms $a_{j,k}$ of (1.1) be given in Figure 4 below. 
In this case, (1.1) fails to converge in Pringsheim's sense, while 
$$\sum^M_{j=m} \ \sum^N_{k=n} a_{j,k} = 0 \quad {\rm if} \quad \min\{m, n\} \ge 1\quad 
{\rm and} \quad M\ge m, \ N\ge n.$$
\bigskip

$$
\zmatrix{
\zv & \vdots & \zv & & \zv & & \zv  & & \zv & & \zv & & \zv & \dddots  \cr
\vonal
\zv & -1 & \zv & 0 & \zv & 0 & \zv & 0 & \zv & 0 & \zv & 0 & \zv  & \cr
\vonal
\zv & 1 & \zv & 0 & \zv & 0 & \zv & 0 & \zv & 0 & \zv & 0 & \zv  & \cr
\vonal
\zv & -1 & \zv & 0 & \zv & 0 & \zv & 0 & \zv & 0 & \zv & 0 & \zv  & \cr
\vonal
\zv & 1 & \zv & 0 & \zv & 0 & \zv & 0 & \zv & 0 & \zv & 0 & \zv  & \cr
\vonal
\zv & -1 & \zv & 0 & \zv & 0 & \zv & 0 & \zv & 0 & \zv & 0 & \zv  & \cr
\vonal
\zv & 1 & \zv & -1 & \zv & 1 & \zv & -1 & \zv & 1 & \zv & -1 & \zv & \cdots & \cr
\vonal
}$$
\abra {4} 

\bigskip

Hardy [1] introduced the notion of {\it convergence in the regular sense}, or briefly: 
the notion  of {\it regular convergence}, of double series as follows: (1.1) is said to converge regularly to the sum 
$s\in \C$  if it converges to $s$ in Pringsheim's sense and, in addition, each of its row and column  series 
defined in (1.4) also converges as a single series. This notion 
was rediscovered by the present author in [4], where it was defined by an equivalent condition, and it was 
called there as `{\it convergence in a  restricted sense}'. Our definition reads as follows: (1.1) 
is said to converge in a restricted sense if for every $\e>0$ there exists 
$\kappa_2 = \kappa_2 (\e) \in \N$ such that 
$$\Big|\sum^M_{j=m} \ \sum^N_{k=n} a_{j,k} \Big| < \e\quad {\rm if} \quad \max\{m,n\} > \kappa_2 
\quad {\rm and} \quad M\ge m\ge 0, \quad N\ge n\ge 0.\leqno(1.6)$$

The essential difference between (1.5) and (1.6) is that the term `min' is changed for `max' relating to the 
lower limits of the double sums involved in them.

It follows from (1.6) that 
$$|s_{m_1, n_1} - s_{m_2, n_2}| < 2\e\quad {\rm if} \quad \min\{m_1, n_1, m_2, n_2\} > \kappa_2(\e).$$
Since this is true for every $\e>0$, by the Cauchy convergence criterion, the convergence of 
(1.1) in a  restricted sense implies its convergence in Pringsheim's sense, and therefore the sum of a 
 double series which satisfies condition (1.6) is well defined. 
Furthermore, taking $M=m$ or 
$N=n$ in (1.6), we see that each of the row and column series in (1.4) satisfies the ordinary 
Cauchy convergence criterion; consequently, each of them converges. This shows that our definition 
(1.6) implies the fulfillment of the definition of regular convergence given by Hardy. To check the 
converse implication, it is enough to  refer to condition (1.5) and to the fact that a finite 
number of convergent row and column series converge uniformly. From now on, we will use the term `regular convergence' 
for double series (1.1) which satisfy condition (1.6). 

It is obvious that if the double series (1.1) {\it converges absolutely}, that is, if 
$$\sum^\infty_{j=0} \ \sum^\infty_{k=0} |a_{j,k}| < \infty,$$
then it also converges regularly. The converse implication is not true in general. See, e.g., the following 

\noindent {\bf Example 5.}  Let the terms $a_{j,k}$ of (1.1) be given in Figure 5. In this case (1.1) converges regularly, 
but it fails to converge absolutely. 

\bigskip

$$
\zmatrix{
\zv & \vdots & \zv & & \zv & & \zv  & & \zv & & \zv & & \zv & \dddots  \cr
\vonal
\zv & 0 & \zv & 0 & \zv & 0 & \zv & 0 & \zv & 0 & \zv & {1\over 6} & \zv  & \cr
\vonal
\zv & 0 & \zv & 0 & \zv & 0 & \zv & 0 & \zv & {1\over 5} & \zv & -{1\over 5} & \zv  & \cr
\vonal
\zv & 0 & \zv & 0 & \zv & 0 & \zv & {1\over 4} & \zv & -{1\over 4}  & \zv & 0 & \zv  & \cr
\vonal
\zv & 0 & \zv & 0 & \zv & {1\over 3} & \zv & -{1\over 3} & \zv & 0 & \zv & 0 & \zv  & \cr
\vonal
\zv & 0 & \zv & {1\over 2} & \zv & - {1\over 2} & \zv & 0 & \zv & 0 & \zv & 0 & \zv  & \cr
\vonal
\zv & 1 & \zv & -1 & \zv & 0 & \zv & 0 & \zv & 0 & \zv & 0 & \zv & \cdots & \cr
\vonal
}$$
\abra {5}

It is clear that if either $\{a_{j,k}\} \subset \overline{\R}_+$ or $\{a_{j,k}\}\subset \overline{\R}_-:= 
(-\infty, 0]$, then the convergence of the double series (1.1) in Pringsheim's sense is equivalent with its absolute convergence, and a fortiori, with its regular 
convergence. 

We recall that if the double series (1.1) converges absolutely, then its sum can be also  computed by successive 
summation as follows: 
$$\sum^\infty_{j=0} \ \sum^\infty_{k=0} a_{j,k} = \sum^\infty_{j=0} \Big(\sum^\infty_{k=0} a_{j,k}\Big) 
= \sum^\infty_{k=0} \Big(\sum^\infty_{j=0} a_{j,k}\Big).\leqno(1.7)$$

The following Theorem 1 is  folklore. 

\bigskip

\noindent {\bf Theorem 1.} {\it If the double series} (1.1) {\it converges regularly, then} (1.7) {\it holds true.} 

For the reader's convenience, in the sequel we present a {\it proof of Theorem} 1. Since the double series (1.1) converges regularly,  it converges in 
Pringsheim's sense, and each of its row and column series also converges. 
Denote by $s$ the sum of the double series 
(1.1), and by $r_j$ the sum of its $j$th row  series, that is, 
$$\sum^\infty_{k=0} a_{j,k} = r_j, \quad j\in \N.\leqno(1.8)$$
Given any $\e>0$, by  (1.6), we have 
$$\Big|\sum^m_{j=0} \ \sum^N_{k=n} a_{j,k}\Big|< \e \quad {\rm if} \quad m\in \N \quad {\rm and} \quad 
N> n > \kappa_2 (\e).$$ 
Letting $N\to \infty$ gives 
$$\Big|\sum^m_{j=0} \ \sum^\infty_{k=n+1} a_{j,k} \Big| \le \e \quad {\rm if} \quad m\in \N 
\quad {\rm and} \quad n > \kappa_2.\leqno(1.9)$$
By (1.8), we may write that 
$$\sum^\infty_{k=n+1} a_{j,k} = r_j - \sum^n_{k=0} a_{j,k}.\leqno(1.10)$$
It follows from (1.2), (1.9) and (1.10) that 
$$\Big|\sum^m_{j=0} r_j - s_{m,n}\Big| = \Big|\sum^m_{j=0} r_j - \sum^m_{j=0} \ \sum^n_{k=0} a_{j,k}\Big| 
\le \e \quad {\rm if} \quad m\in \N \quad {\rm and} \quad n> \kappa_2.$$
Combining the inequality just obtained with (1.3) yields 
$$\Big|\sum^m_{j=0} r_j - s\Big| \le \Big|\sum^m_{j=0} r_j - s_{m,n}\Big| + |s_{m,n} - s| < 2\e$$
$${\rm if} \quad \min\{m,n\} > \max \{\kappa_1, \kappa_2\}.$$
Since $\e>0$ is arbitrary, this proves the first equality in (1.7). 

The proof of the second equality in (1.7) runs along analogous lines. The proof of (1.7) is complete. $\kocka$ 
\bigskip

\noindent  {\bf Remark 1.} The present author in a joint paper [2] with P. K\'orus studied the uniform convergence of the family of double sine integrals of the form 
$$\sum^\infty_{j=1} \ \sum^\infty_{k=1} a_{j,k} \sin jx \sin ky, \quad {\rm where} \quad 
(x,y) \in [0, 2\pi)^2\leqno(1.11)$$
and the double sequence $\{a_{j,k}\} \subset \C$ satisfies one of  certain generalized monotonicity conditions. 
Now, we proved the following 

\bigskip

\noindent {\bf Theorem A} (see [2, Theorem 1]). {\it  If $\{a_{j,k}\} \subset \C$ is such that } 
$$jk a_{j,k} \to 0 \quad as\quad \max \{j,k\}\to \infty, \leqno(1.12)$$
{\it then the regular convergence of the double sine series} (1.11) {\it is uniform in} $(x,y)$. 

{\it Conversely, if $\{a_{j,k}\} \subset \overline{\R}_+$ is such that the regular convergence of the double sine series} 
(1.11) {\it is uniform in} $(x,y)$, {\it then condition} (1.12) {\it is satisfied.} 

On the other hand, all our attempts have failed so far to prove an analogue of Theorem A when 
regular convergence is exchanged  for convergence in Pringsheim's sense, and the `max' is exchanged for `min' in (1.12).  
To be more precise, the proof of the sufficiency part of any analogue of Theorem A seems to be hopeless in the case of convergence 
in Pringsheim's sense uniformly in $(x,y)$; while the proof of the necessity part of the analogue of Theorem A is 
routine. 
\bigskip

\noindent {\bf Remark 2.} In Harmonic Analysis (e.g., Fourier series, see in [9, Ch. XVII] by A. Zygmund) 
one frequently meets double series of the form 
$$\sum^\infty_{j=-\infty} \ \sum^\infty_{k=-\infty} a_{j,k}, \leqno(1.13)$$
where $\{a_{j,k}: (j,k) \in \Z^2\} \subset \C$. Using the notion of  {\it symmetric rectangular partial sum} 
defined by 
$$s_{m,n} : = \sum^m_{j=-m} \ \sum^n_{k=-n} a_{j,k} : = \sum_{|j|\le m} \ 
\sum_{|k|\le n} a_{j,k}, \quad (m,n) \in \N^2,$$
the convergence of (1.13) in Pringsheim's sense is also defined by (1.3). Accordingly, the so-called row series 
of (1.13) are defined by 
$$\sum^\infty_{k=-\infty} a_{0, k} \quad {\rm and} \quad \sum^\infty_{k=-\infty} (a_{j,k} + a_{-j,k}), 
\quad j\in \N_+;$$
and the column series of (1.13) are defined analogously (cf. (1.4)). 

In the definition of the regular convergence of (1.13), instead of (1.6) we require the fulfillment of the following condition: 
$$\Big|\sum_{m\le |j| \le M} \ \sum_{n\le |k|\le N} a_{j,k}\Big| < \e \quad {\rm if} \quad 
\max\{m,n\} > \kappa_2 (\e) \quad {\rm and} \quad M\ge m\ge 0, N\ge n\ge 0.$$

\bigskip

\noindent{\bf Remark 3.} Similarly to the above two convergence notions for double series, one can 
define analogous convergence notions for double sequences of complex numbers. The double sequence 
$(s_{m,n}) = (s_{m,n}: (m,n) \in \N^2) \subset \C$ is said to {\it converge in Pringsheim's sense} to some 
$s\in \C$, in symbols:
$$\lim_{m,n\to \infty} s_{m,n} = s,$$
if for every $\e>0$ there exists some $\kappa_1 = \kappa_1(\e) \in \N$ such that condition 
(1.3) is satisfied. 

Next, taking into account the notation (1.2), we may write that 
$$\sum^M_{j=m} \ \sum^N_{k=n} a_{j,k} = s_{M,N} - s_{m-1, N} - s_{M, n-1} + s_{m-1, n-1}$$
with the agreement that 
$$s_{m,n} := 0 \quad {\rm if}\quad \min\{m,n\} = -1, \quad (m+1, n+1) \in \N^2.$$
Keeping definition (1.6) in mind , the double sequence $(s_{m,n})$ is said to {\it converge regularly} if for every $\e>0$ there exists some 
$\kappa_2 = \kappa_2 (\e) \in \N$ such that 
$$|s_{M,N} - s_{m-1, N} - s_{M, n-1} + s_{m-1, n-1}| < \e$$
$${\rm if} \quad \max\{m,n\} > \kappa_2 \quad {\rm and} \quad M\ge m\ge 0, \quad N\ge n\ge 0.$$

It is routine to check that the regular convergence of a double sequence implies its convergence in 
Pringsheim's sense. Therefore, the finite limit of a regularly convergent double sequence is well defined. 
Furthermore, if a double sequence $(s_{m,n})$ converges regularly, then each of the single sequences 
$(s_{m,n} : m\in \N)$ for fixed $n\in \N$, and $(s_{m,n} : n\in \N)$ for fixed $m\in \N$, 
converges; and similarly to Theorem 1, we have 
$$\lim_{m\to \infty}( \lim_{n\to \infty}s_{m,n}) = \lim_{n\to \infty} (\lim_{m\to \infty} s_{m,n}) = \lim_{m,n\to \infty} 
s_{m,n} =: s,$$
where $s$ is the finite limit of the double sequence $(s_{m,n})$ in Pringsheim's sense.

\vglue1cm

\noindent {\bf 2. Convergence of double integrals of locally integrable functions over $\overline{\R}^2_+$} 

We  consider the double integral 
$$\int^\infty_0 \int^\infty_0 f(u,v) du dv\leqno(2.1)$$
over the closed first quadrant $\overline{\R}^2_+ := [0, \infty)^2$. The  {\it rectangular partial integrals} of (2.1) are defined by 
$$I(x,y): = \int^x_0 \int^y_0 f(u,v) du dv, \quad (x,y) \in \R^2_+,\leqno(2.2)$$
where $f: \overline{\R}^2_+ \to \C$ is a locally integrable function in Lebesgue's sense; that is,  
$f\in L^1 ({\cal R})$ for every bounded rectangle 
$${\cal R} := [a_1, b_1] \times [a_2, b_2], \quad 0\le a_\ell < b_\ell < \infty, \ \ \ell=1,2;\leqno(2.3)$$
in symbols: $f\in L^1_{\loc} (\overline{\R}^2_+)$. 

Analogously to the case of double series, the double integral (2.1) is said to {\it converge in Pringsheim's sense} 
to some  $I\in \C$ , 
or it is equivalently said that $I$ is {\it the sum of} (2.1) {\it in Pringsheim's sense}, in symbols:  
$$\lim_{x,y\to \infty} I(x,y) = I, \leqno(2.4)$$
if for every $\e>0$ there exists
$\rho_1 = \rho_1(\e) \in \R_+$ such that 
$$|I(x,y) - I| < \e\quad {\rm if} \quad \min\{x,y\} > \rho_1. \leqno(2.5)$$

We note that if the double integral (2.1) converges in Pringsheim's sense, then it follows from (2.2) and (2.5) that for every $\e>0$, 

$$\Big|\int^{x_1}_x \int^{y_1}_y f(u,v) du dv\Big| = |I(x_1, y_1) - I(x,y_1) - I(x_1, y) + I(x,y)| < 4\e \leqno(2.6)$$
$${\rm if} \quad \min\{x,y\} > \rho_1 \quad {\rm and} \quad x_1 > x\ge 0, \ y_1 > y\ge 0.$$
Clearly, (2.6) is the nondiscrete counterpart of (1.5). 
The converse statement is not true in general.  For example, consider the double series (1.1) in Example 4 of Section 1, and define 
$$f(u, v) : = a_{j,k} \quad {\rm if} \quad (u,v) \in [j, j+1) \times [k, k+1), \quad (j,k) \in \N^2.\leqno(2.7)$$
Then the double integral (2.1) fails to converge in Pringsheim's sense, while (2.6) is satisfied since 
$$\int^{x_1}_x \int^{y_1} _y f(u,v) du dv = 0 \quad {\rm if} \quad \min \{x,y\} \ge 1 \quad {\rm and} 
\quad x_1 > x\ge 0, \ y_1 > y\ge 0.$$

We recall (see in [5]) that the double integral (2.1) is said to 
{\it converge regularly} if for every $\e>0$ there exists $\rho_2 = \rho_2(\e) \in \R_+$ such that 
$$\Big|\int^{x_1}_x \int^{y_1}_y f(u,v) du dv\Big| < \e \quad {\rm if} \quad \max \{x,y\} > \rho_2\leqno(2.8)$$ 
$${\rm and} \quad x_1 > x\ge 0, \quad y_1> y \ge 0.$$
The essential difference between (2.5) and (2.8) is that 
the term `min' is changed for the term `max' relating to the lower limits in the double integrals involved in 
them. 

It follows from (2.8) that 
$$|I(x_1, y_1) - I(x_2, y_2)| < 2\e \quad {\rm if} \quad \min \{x_1, y_1, x_2, y_2\} > \rho_2 (\e).$$
Since this is true for every $\e>0$, by the Cauchy convergence criterion, the regular convergence of (2.1) implies its convergence in 
Pringsheim's sense, and therefore the sum of a regularly convergent double integral (2.l) is well defined. 
The converse statement is not true in general. For example, consider the double series (1.1) in 
Example 1 of Section 1, and define the function $f$ by (2.7). Then  the double integral (2.1) converges to 
0 in Pringsheim's sense, but it fails to converge regularly. 

It is obvious that if $f\in L^1 (\overline{\R}^2_+)$, then the double integral (2.1) converges regularly, 
and its sum $I$ equals the integral of $f$ over the whole first quadrant $\overline{\R}^2_+$. However, the double integral (2.1) may 
converge regularly in such cases when 
$f\not\in L^1 (\overline{\R}^2_+)$, as the following example shows.  

\bigskip

\noindent {\bf Example 6}. Consider  the terms $a_{j,k}$ of the double series (1.1) given in Figure 6  below, 
and define the function $f$ by (2.7). Then  the double integral (2.1) converges regularly, although 
$f\not\in L^1 (\overline{\R}^2_+)$. In this case, even the marginal functions 
$f(\cdot, v_0) \not\in L^1 (\overline{\R}_+)$ and $f(u_0, \cdot) \not \in L^1 (\overline{\R}_+)$ for any fixed  $v_0\in \overline{\R}_+$ 
and $u_0 \in \overline{\R}_+$, respectively. 

\bigskip

$$
\zmatrix{
\zv & \vdots & \zv & & \zv & & \zv  & & \zv & & \zv & & \zv & \dddots  \cr
\vonal
\zv & -{1\over 3}  & \zv & {1\over 3}  & \zv & -{1\over 4}  & \zv & {1\over 4}  & \zv & -{1\over 5}  & \zv & {1\over 5} & \zv  & \cr
\vonal
\zv & {1\over 3}  & \zv & -{1\over 3}  & \zv & {1\over 4}  & \zv & -{1\over 4} & \zv & {1\over 5} & \zv & -{1\over 5} & \zv  & \cr
\vonal
\zv & -{1\over 2} & \zv & {1\over 2}  & \zv & -{1\over 3} & \zv & {1\over 3} & \zv & -{1\over 4} & \zv & {1\over 4} & \zv  & \cr
\vonal
\zv & {1\over 2}  & \zv & -{1\over 2} & \zv & {1\over 3} & \zv & -{1\over 3}  & \zv & {1\over 4}  & \zv & -{1\over 4} & \zv  & \cr
\vonal
\zv & -1 & \zv & 1 & \zv & -{1\over 2}  & \zv & {1\over 2}  & \zv & -{1\over 3}  & \zv & {1\over 3}  & \zv  & \cr
\vonal
\zv & 1 & \zv & -1 & \zv & {1\over 2}  & \zv & -{1\over 2}  & \zv & {1\over 3}  & \zv & -{1\over 3}  & \zv & \cdots & \cr
\vonal
}$$
\abra {6}

\bigskip

\noindent {\bf Remark 4.} If the double integral (2.1) converges regularly, then the finite 
limit of the so-called {\it `horizontal strip' integral} 
$$\int^{x_1}_0 \int^{y_1}_y f(u,v) du dv \quad {\rm exists\ as} \quad x_1\to \infty, \leqno(2.9)$$
locally uniformly in $(y, y_1)$, where $y_1 > y\ge 0$. By the term `locally uniform convergence' in this case, we mean that for every 
$c\in \R_+$ the finite limits in (2.9) exist uniformly in $(y, y_1)$, where $0\le y < y_1\le c$. In other words, for every 
$c>0$ and $\e>0$, there exists $\rho_3 = \rho_3 (c, \e) > 0$ such that 
$$\Big|\int^{x_1}_0 \int^{y_1}_y \ f(u,v) du dv\Big| < \e \quad {\rm for\ all} \quad 0\le y < y_1 \le c 
\quad {\rm if} \quad x_1 > \rho_3.$$

Analogously, if the double integral (2.1) converges regularly, then the finite limit of the so-called `{\it vertical strip' integral} 
$$\int^{x_1}_x \int^{y_1}_0 f(u,v) du dv \quad {\rm exists\ as} \quad y_1\to \infty, \leqno(2.10)$$
locally uniformly in $(x, x_1)$, where $x_1> x\ge 0$. 

It is worth formulating the following characterization of regular convergence of double integrals, whose proof is routine. 

\bigskip

\noindent {\bf Theorem 2.} {\it Suppose $f\in  L^1_{\loc} (\overline{\R}^2_+)$. The double integral} (2.1) 
{\it converges regularly if and only if} 

(i) {\it it converges in Pringsheim's sense;} 

(ii) {\it the finite limit of the `horizontal strip' integral in} (2.9) {\it exists locally uniformly in $(y, y_1)$, where} 
$y_1 > y \ge 0$; 
 {\it as well as the finite limit of the `vertical strip' integral in} (2.10) 
{\it exists locally uniformly in} $(x, x_1)$, {\it where} $x_1> x \ge 0$. 

\noindent {\it Proof of Theorem} 2. It is routine. 

In the following Example 7 we show that (contrary to the case of double series) the sufficiency part of Theorem 2 fails if we drop the requirement 
``locally uniformly" in condition (ii). 

\bigskip

\noindent {\bf Example 7.} We define the function $f(u,v)$ on the lower 
half $\{(u, v): 0 \le v \le u\}$ of $\overline{\R}^2_+$ as follows
$$f(u,v):= \cases{1 \quad &if\quad $(u,v) \in [2^k, 3\cdot 2^{k-1}) \times (3\cdot 2^{-k-2}, 2^{-k}]$\cr
&{} \quad $\bigcup [3\cdot 2^{k-1}, 2^{k+1}) \times (2^{-k-1}, 3\cdot 2^{-k-2}]$,\cr
-1 &if\quad $(u,v) \in [3\cdot 2^{k-1}, 2^{k+1}) \times (3\cdot 2^{-k-2}, 2^{-k}]$\cr
&{} \quad $ \bigcup [2^k, 3\cdot 2^{k-1}) 
\times (2^{-k-1}, 3\cdot 2^{-k-2}]$, \ {\rm where}\ $k=0,1,2,\ldots$;\cr
0 &otherwise;\cr}$$

\hskip 1,3cm while on the upper half $\{(u,v): 0\le v < u\}$ of $\overline{\R}^2_+$ we set $f(u,v):= f(v,u)$. 

This $f\in L^1_{\loc} (\overline{\R}^2_+)$, the double integral (2.1) converges to 
$0$ in Pringsheim's sense, all its 
`horizontal strip' integrals as well as all its `vertical strip' integrals converge to 0, but not locally 
uniformly, and the double integral (2.1) fails to onverge regularly. We even observe that the marginal 
functions $f(\cdot, v_0) \in L^1(\overline{\R}_+)$ and $f(u_0, \cdot) \in L^1(\overline{\R}_+)$ for any fixed 
$v_0 \in \overline{\R}_+$ and $u_0 \in \overline{\R}_+$, respectively; and their integrals over 
$\overline{\R}_+$ are equal to 0. 

It is clear that if the real-valued function $f\in L^1_{\loc} (\overline{\R}^2_+)$ is such that $f(u,v) \ge 0$ almost 
everywhere on $\overline{\R}^2_+$, then from the convergence of the double integral (2.1) in Pringsheim's sense it follows that 
$f\in L^1 (\overline{\R}^2_+)$; thus, in this case the convergence of (2.1) in Pringsheim's sense is equivalent with its regular 
convergence. 

\vfill\eject

\noindent {\bf Remark 5.} The present author in a joint paper [3] with P. K\'orus studied the uniform convergence of the family of {\it double 
sine integrals} of the form 
$$\int^\infty_0 \int^\infty_0 f(u,v) \sin ux \sin vy du dv, \quad {\rm where} \quad 
(x,y) \in \R^2_+, \leqno(2.11)$$
$f: \R^2_+ \to \C$ is a measurable function in Lebesgue's sense, which 
satisfies one of certain generalized monotonicity 
conditions and the condition that 
$$u v f(u,v) \in L^1 (\overline{\R}^2_+).$$
Now, we proved the following 

\bigskip

\noindent {\bf Theorem B} (see [3, Theorems 1 and 2]). {\it If $f: \R^2_+ \to \C$ is such that 
$$uv f(u,v) \to 0 \quad as\quad \max\{u,v\} \to \infty \leqno(2.12)$$
and
$${1\over xy} \int^x_0 \int^y_0 uv f(u,v) du dv \to 0 \quad as\quad \max\{x,y\} \to \infty,\leqno(2.13)$$
then the regular convergence of the double sine integrals} (2.11) {\it is uniform in $(x,y)$. 

Conversely, if $f: \R^2_+ \to \overline{\R}_+$ is such that $uv f(u,v) \in L^1_{\loc} (\overline{\R}^2_+)$ and  the regular 
convergence of the double sine integrals} (2.11) {\it is uniform in $(x,y)$, then condition } (2.12) 
{\it is satisfied. }

On the other hand, all our attemps have failed so far to prove an analogue of Theorem B when 
regular convergence is exchanged  for convergence in Pringsheim's sense, 
and the `max' is exchanged for `min' in (2.12) and (2.13). 
Similarly to our remark made after Theorem A, the proof of the sufficiency part of any analogue of 
Theorem B seems to be hopeless in the case of convergence in Pringsheim's sense uniformly in $(x,y)$; 
while the proof of the necessity part of the analogue of Theorem B is routine. 

\bigskip

\noindent {\bf Remark 6.} In Harmonic Analysis (e.g., Fourier transform, see in [8, Ch. I] by E. M. Stein and G. Weiss) 
one frequently meets double integrals of the form 
$$\int^\infty_{-\infty} \int^\infty_{-\infty} f(u,v) du dv, \quad {\rm where} \quad f\in L^1_{\loc} (\R^2).\leqno(2.14)$$
Using the notion of {\it symmetric rectangular partial integral} defined by 
$$I(x,y) : = \int^x_{-x} \int^y_{-y} f(u,v) du dv, \quad (x,y) \in \R^2_+,$$
the convergence of (2.14) in Pringsheim's sense is also defined by (2.5). In the definition of the regular convergence of 
(2.14), instead of (2.8) we require the fulfillment of the following condition: for every $\e>0$ there 
exists $\rho_2 = \rho_2 (\e) \in \R_+$ such that 
$$\Big|\int_{x< |u| < x_1} \int_{y<|v|<y_1} f(u,v) du dv\Big| < \e \quad {\rm if} \quad \max\{x,y\} > \rho_2$$
$${\rm and} \quad x_1 > x\ge 0, \quad y_1 > y\ge 0.$$

\vglue1cm

\noindent {\bf 3. Main result: a generalized version of Fubini's theorem} 

By Fubini's theorem on successive integration (see, e.g., in [7, on p. 85]), if 
$f\in L^1 (\overline{\R}^2_+)$, then the marginal functions $f(\cdot, v_0)$ and $f(u_0, \cdot)$ belong to the space 
$L^1 (\overline{\R}_+)$ for almost all fixed $v_0 \in \R_+$ and $u_0 \in \R_+$, respectively, and 
$$\int^\infty_0 \int^\infty_0 f(u,v) du dv = \int^\infty_0 \Big(\int^\infty_0 f(u,v) dv\Big) du = 
\int^\infty_0 \Big(\int^\infty_0 f(u,v) du\Big) dv. \leqno(3.1)$$

Before formulating our main Theorem 3 below, we make the following observation. If $f\in L^1_{\loc} (\overline{\R}^2_+)$, then 
$f\in L^1 ({\cal R})$ for every bounded rectangle of form (2.3). By Fubini's theorem, the marginal function 
$f(\cdot, v_0) \in L^1 ([a_1, b_1])$ for almost all fixed $v_0 \in [a_2, b_2]$, and  the marginal function 
$f(u_0, \cdot) \in L^1 ([a_2, b_2])$ for almost all fixed  
$u_0 \in [a_1, b_1]$. Since a countable union of sets of measure zero is also of measure zero, it follows that if $f\in L^1_{\loc} 
(\overline{\R}^2_+)$, then the marginal functions $f(\cdot, v_0)$ and $f(u_0, \cdot)$ belong to the class $L^1_{\loc} (\overline{\R}_+)$ for 
almost all fixed $v_0 \in \R_+$ and $u_0 \in \R_+$, respectively. 

Now, our main result reads as follows. 

\bigskip

\noindent {\bf Theorem 3.} {\it If $f\in L^1_{\loc} (\overline{\R}^2_+)$ and its double integral} (2.1) {\it converges regularly, then the finite 
limits
$$\lim_{y\to \infty} \int^A_0 \Big(\int^y_0 f(u,v) dv\Big) du =: I_1 (A), \quad 0<A<\infty, \leqno(3.2)$$
and
$$\lim_{x\to \infty} \int^B_0 \Big(\int^x_0 f(u,v) du \Big) dv = : I_2 (B), \quad 0<B<\infty, \leqno(3.3)$$
exist uniformly in $A$ and B, respectively;  and we have 
$$\lim_{A\to \infty} I_1(A) = \lim_{B\to \infty} I_2 (B) = 
\lim_{A,B\to \infty} \int^A_0 \int^B_0 f(u,v) du dv=: I,\leqno(3.4)$$
where $I\in \C$ is the sum of the double integral} (2.1) {\it in Pringsheim's sense}. 

In order to reformulate Theorem 3 in a transparent form, first we introduce the following formal notations: 
$$\int^{\to \infty}_0 \int^{\to \infty}_0 f(u,v) du dv:= I,$$
$$\int^A_0 \Big(\int^{\to \infty}_0  f(u,v) dv\Big) du: = I_1 (A),$$
$$\int^B_0 \Big(\int^{\to \infty}_0  f(u,v) du\Big) dv: = I_2 (B),$$
$$\int^{\to \infty}_0 \Big(\int^{\to \infty}_0 f(u,v) dv\Big) du := \lim_{A\to \infty} I_1 (A),$$
$$\int^{\to \infty}_0  \Big(\int^{\to \infty}_0  f(u,v) du\Big) dv := \lim_{B\to \infty} I_2(B), \quad 
{\rm where} \quad 0 < A,B < \infty.$$ 

Using these notations, (3.4) can be rewritten into the following form: 
$$\eqalign{\int^{\to \infty}_0 \int^{\to \infty}_0 f(u,v) du dv &= \int^{\to \infty}_0 \Big(\int^{\to \infty}_0 f(u,v) dv\Big) du\cr
&= \int^{\to \infty}_0 \Big(\int^{\to \infty}_0 f(u,v) du\Big) dv.\cr}\leqno(3.5)$$
Clearly, (3.5) can be considered as a generalized version of Fubini's theorem on 
successive integration of the double integral (2.1), provided that 
$f\in L^1_{\loc} (\overline{\R}^2_+)$ and that the double integral (2.1) converges regularly. 

The following corollary of Theorem 3 is of  interest. 

\bigskip

\noindent {\bf Corollary 4.} {\it Suppose $f\in L^1_{\loc} (\overline{\R}^2_+)$ and the marginal functions 
$f(\cdot, v_0)$ and $f(u_0, \cdot)$ belong to the space $L^1 (\overline{\R}_+)$ for almost all fixed 
$v_0 \in \R_+$ and $u_0 \in \R_+$, respectively. If the double integral} (2.1) {\it converges regularly, then 
$$I_1 (A) = \int^A_0 \Big(\int^\infty_0 f(u,v) dv\Big) du, \ \ I_2 (B) = \int^B_0 \Big(\int^\infty_0 f(u,v) du\Big) dv, $$
and} (3.4) {\it holds true.} 

We note that under the conditions of Corollary 4, (3.5) can be rewritten in the following form: 
$$\eqalign{\int^{\to \infty}_0 \int^{\to \infty}_0 f(u,v) du dv &= \int^{\to \infty}_0\Big(\int^\infty_0 
f(u,v) dv\Big) du\cr
&= \int^{\to \infty}_0\Big(\int^\infty_0 f(u,v) du\Big) dv.\cr}\leqno(3.6)$$

If $f\in L^1 (\overline{\R}^2_+)$, then (3.6) clearly coincides with (3.1). 

\noindent {\it Proof of Theorem} 3. By the definition (2.8) of regular convergence, for every $\e>0$ 
there exists $\rho_2 = \rho_2 (\e) \in \R_+$ such that 
$$\Big|\int^A_0 \int^y_0 f(u,v) du dv - \int^A_0 \int^{y_1}_0 f(u,v) du dv\Big|\leqno(3.7)$$
$$=\Big|\int^A_0 \int^y_{y_1} f(u,v) du dv\Big| < \e\quad {\rm if} \quad A> 0 \quad {\rm and} \quad y> y_1 > \rho_2(\e).$$
Applying Fubini's theorem on the successive integration over the bounded rectangle 
${\cal R} : = [0, A] \times [y_1, y]$, inequality (3.7) can be rewritten into the following form: 
$$\sup_{A>0} \Big|\int^A_0 \Big(\int^y_0 f(u,v) dv\Big) du\leqno(3.8)$$
$$-\int^A_0 \Big(\int^{y_1}_0 f(u,v) dv\Big) du \Big| \le \e \quad {\rm if} \quad y> y_1 > \rho_2 (\e).$$
Since $\e>0$ is arbitrary in (3.8), by virtue of the Cauchy convergence criterion, the finite limit in (3.2) exists, which 
we denote by $I_1 (A)$; and this limit exists even uniformly in $A>0$. 

An analogous argument proves the existence of the finite limit in (3.3), even uniformly in $B>0$ . 

Letting $y\to \infty$ in (3.8) gives 
$$\sup_{A>0} \Big|I_1 (A) - \int^A_0 \Big(\int^{y_1}_0 f(u,v) dv \Big) du \Big| \le \e \quad {\rm if} \quad 
y_1 > \rho_2 (\e).$$
Setting $y_1:= B$ in the inequality just obtained, we find that  for all $0<A<\infty$, 
$$\Big|I_1 (A) - \int^A_0 \Big(\int^B_0 f(u,v) dv\Big) du \Big| \le \e \quad 
{\rm if} \quad B> \rho_2 (\e).\leqno(3.9)$$

The symmetric counterpart of (3.9) can be proved in an analogous way: for all $0<B<\infty$, 
$$\Big|I_2 (B) - \int^B_0 \Big(\int^A_0 f(u,v) du \Big) dv \Big| \le \e \quad {\rm if} \quad A> \rho_2 (\e).\leqno(3.10)$$

Since $f\in L^1_{\loc} (\overline{\R}^2_+)$, we may apply Fubini's theorem on successive integration over the bounded rectangle 
$[0, A] \times [0, B]$ to obtain that 
$$\int^A_0 \int^B_0 f(u,v) du dv = \int^A_0 \Big(\int^B_0 f(u,v) dv\Big) du\leqno(3.11)$$
$$= \int^B_0 \Big(\int^A_0 f(u,v) du\Big) dv \quad {\rm for\ all}\quad 0<A,B<\infty.$$
In particular, it follows from (2.5) and (3.9) - (3.11) that 
$$|I_1 (A) - I| \le |I_1 (A) - \int^A_0 \Big(\int^B_0 f(u, v) dv\Big) du|\leqno(3.12)$$
$$ + 
\Big|\int^A_0 \int^B_0 f(u,v) du dv - I\Big| < 2\e,$$
and analogously we find that 
$$|I_2(B) - I| < 2\e 
 \quad {\rm if} \quad 
\min\{A,B\} > \max\{\rho_1 (\e), \rho_2 (\e)\}.\leqno(3.13)$$
Since $\e>0$ is arbitrary in (3.12) and (3.13),  we conclude (3.4) to be proved. $\kocka$ 
\bigskip

{\bf Acknowledgements.} The author is indebted to Professor John Horv\'ath for careful reading of the manuscript 
and for valuable suggestions to improve the presentation.

\vglue1cm

{\bf References} 

\item{[1]} G. H. Hardy, On the convergence of certain multiple series, 
{\it Proc. Cambridge Philosoph. Soc.}, {\bf 19} (1916-1919), 86-95. 

\item{[2]} P. K\'orus and F. M\'oricz, On the  uniform convergence of double sine series, 
{\it Studia Math.}, {\bf 193} (2009), 79-97. 

\item{[3]} P. K\'orus and F. M\'oricz, Generalizations to monotonicity for uniform convergence of double sine integrals over 
$\overline{\R}^2_+$, {\it Studia Math.}, {\bf 201} (2010), 287-304. 

\item{[4]} F. M\'oricz, On the convergence in a restricted sense of multiple series, 
{\it Analysis Math.}, {\bf 5} (1979), 135-147. 

\item{[5]} F. M\'oricz, On the uniform convergence of double sine integrals over $\overline{\R}^2_+$, 
{\it Analysis}, {\bf 31} (2011), 191-204.

\item{[6]} A. Pringsheim, Zur Theorie der zweifach unendlichen Zahlenfolgen, Math. Ann. 53(1900), 289-321.

\item{[7]} F. Riesz et B. Sz.-Nagy, {\it Lecons d'analyse fonctionelle}, Gauthier-Villars, 
Paris, 1955. 

\item{[8]} E. M. Stein and G. Weiss, {\it Introduction to Fourier Analysis on Euclidean Spaces}, 
Princeton Univ. Press, 1971. 

\item{[9]} A. Zygmund, {\it Trigonometric series}, Vol. II, Cambridge Univ. Press, 1959. 

\vglue1cm

\noindent updated version, December 19, 2011.

\bye